\documentclass[a4paper,10pt]{article}
\usepackage{amsmath}
\usepackage{amssymb,amsfonts,amsthm}
\usepackage[applemac]{inputenc}
\usepackage[english]{babel}
\usepackage[pdftex]{graphicx}
\usepackage{caption}
\usepackage{subcaption}
\usepackage{anysize}
\usepackage{enumerate}
\usepackage[all]{xy}
\usepackage{mathrsfs}
\usepackage{verbatim}
\usepackage{cite,lpic}

\definecolor{green}{rgb}{0,0.8,0.5}

\usepackage[pdftex]{hyperref}
\hypersetup{colorlinks=true,linkcolor=blue,citecolor=blue}
\usepackage[width=0.8\textwidth]{caption}

\renewenvironment{abstract}{\small\quotation\noindent
 {\bfseries \abstractname .}}{\endquotation \par}

\newenvironment{trivlistparm}[2]{\trivlist
\item[{#1 #2}]}{\endtrivlist}

\newenvironment{prooftext}[1]{\trivlistparm{\bfseries}{#1}}{\Qed\endtrivlistparm}
\newenvironment{prova}{\trivlistparm{\bfseries}{Proof.}}{\Qed\endtrivlistparm}

\catcode`\@=11

\def\resetthefootnote{\renewcommand{\thefootnote}{\@arabic\c@footnote} }
\def\@principiremex#1{\trivlist
 \item[\hskip \labelsep{\bfseries #1\ \thethm.}]\ignorespaces}
\def\opar@principiremex#1[#2]{\trivlist
 \item[\hskip \labelsep{\bfseries #1\ \thethm\ (#2).}]\ignorespaces}

\newcommand{\newTHEOremrom}[2]{\newenvironment{#1}{\refstepcounter{thm}\@ifnextchar[{\opar@principiremex{#2}}
{\@principiremex{#2}}}{\qedB\endtrivlist}} \catcode`\@=12
\DeclareMathSymbol{\square}{\mathord}{AMSa}{"03}
\newcommand{\qedB}{\nopagebreak\hspace*{\fill}$\square$\par}
\newcommand{\Qed}{\nopagebreak\hspace*{\fill}{\vrule width6pt height6pt depth0pt}\par}

\newTHEOremrom{defi} {Definition}
\newTHEOremrom {rem} {Remark}
\newTHEOremrom {ex} {Example}

\newcommand{\refc}[1]{\mbox{$(\ref{#1})$}}
\newcommand{\secc}[1]{Section~\ref{#1}}
\newcommand{\teoc}[1]{Theorem~\ref{#1}}
\newcommand{\propc}[1]{Proposition~\ref{#1}}

\newcommand{\lemc}[1]{Lemma~\ref{#1}}
\newcommand{\defic}[1]{Definition~\ref{#1}}
\newcommand{\obsc}[1]{Remark~\ref{#1}}

\newcommand{\figc}[1]{Figure~\ref{#1}}

\renewcommand{\geq}{\geqslant}
\renewcommand{\epsilon}{\varepsilon}
\renewcommand{\leq}{\leqslant}
\newcommand{\R}{\mathbb{R}}

\newcommand{\Z}{\mathbb{Z}}
\newcommand{\N}{\mathbb{N}}

\newcommand{\sist}[2]{
  \left\{\!
   \begin{array}{l}
    \dot x=#1 \\[2pt] \dot y=#2
   \end{array}
  \right.
}

\newcommand{\PA}{\mathscr{P}}
\newcommand{\PI}{\mathcal{I}}

\newcommand{\RP}{\mathbb{RP}}

\newcommand{\DD}{\mathscr{D}}
\newcommand{\out}{\Pi}
\newcommand{\op}{\ensuremath{\mbox{\rm o}}}

\newtheorem{thm}{Theorem}[section]
\newtheorem{thmx}{Theorem}
\newtheorem{prop}[thm]{Proposition}

\newtheorem{lema}[thm]{Lemma}

\title{\textbf{A criticality result for polycycles in a family  \\Êof quadratic reversible centers}
\footnotetext{Both authors are partially supported by the MINECO Grant MTM2014-52209-C2-1-P. D.~Rojas is also partially supported by the MEC/FEDER grant MTM2014-52232-P.}}

\author{D. Rojas$^{\,1}$ and J. Villadelprat$^{\,2}$\\[10pt]
{\small $^{\,1}$\textsl{Departamento de Matem\'atica Aplicada,}}\\
\vspace{-2pt}
{\small \textsl{Universidad de Granada, 18071 Granada, Spain}}\\[5pt]
{\small $^{\,2}$\textsl{Departament d'Enginyeria Inform\`{a}tica i Matem\`{a}tiques},}\\
\vspace{-2pt}
{\small \textsl{Universitat Rovira i Virgili, Tarragona, Spain}}}

\date{}

\begin{document}

\maketitle

\begin{abstract}
We consider the family of dehomogenized LoudÕs centers $X_{\mu}=y(x-1)\partial_x+(x+Dx^2+Fy^2)\partial_y,$ 
where $\mu=(D,F)\in\R^2,$ and we study the number of critical periodic orbits that emerge or dissapear from the polycycle at the boundary of the period annulus. This number is defined exactly the same way as the well-known notion of cyclicity of a limit periodic set and we call it criticality. The previous results on the issue for the family $\{X_{\mu},\mu\in\R^2\}$ distinguish between parameters with criticality equal to zero (regular parameters) and those with criticality greater than zero (bifurcation parameters). A challenging problem not tackled so far is the computation of the criticality of the bifurcation parameters, which form a set $\Gamma_{B}$ of codimension 1 in $\R^2$. In the present paper we succeed in proving that a subset of $\Gamma_{B}$ has criticality equal to one. 
\end{abstract}

\section{Introduction and statement of the result}\label{section:intro}

In the present paper we study the local bifurcation diagram of the period function associated to a family of quadratic centers. By local we mean near the polycycle at the boundary of the period annulus of the center. In the literature one can find different terminology to classify the quadratic centers but essentially there are four families: Hamiltonian, reversible $Q_3^R$, codimension four $Q_4$ and generalized Lotka-Volterra systems~$Q_3^{LV}$. 
According to Chicone's conjecture~\cite{Chicone}, the reversible centers have at most two critical periodic orbits, whereas the centers of the other three families have monotonic period function. In this context critical periodic orbits play exactly the same role as limit cycles in the celebrated Hilbert's 16th problem (see for instance \cite{Ilyashenko} and references therein). What is more, from the point of view of the techniques, results and notions involved, Chicone's conjecture is the counterpart for the period function to the question of whether quadratic polynomial differential systems have at most four limit cycles, i.e., $H(2)=4.$ Both problems are far from being solved and pose challenging difficulties. There are many papers proving partial results related to Chicone's conjecture and there is much analytic evidence that it is true. In this direction, and without being exhaustive, let us quote Coppel and Gavrilov~\cite{CG}, who showed that the period function of any Hamiltonian quadratic center is monotonic, and Zhao~\cite{Zhao}, who proved the same property for the~$Q_4$ centers. There are very few results concerning the $Q_3^{LV}$ centers. In the middle 80s several
authors \cite{Hsu,Rothe,Waldvogel} showed independently the monotonicity of the classical Lotka-Volterra centers, which constitute a hypersurface inside the $Q_3^{LV}$ family, and more recently the same property has been proved in \cite{Jordi} for two other hypersurfaces. With regard to reversible quadratic centers, it is well known that they can be brought by an affine transformation and a constant rescaling of time to the \emph{Loud normal form}
\begin{equation*}
 \sist{-y+Bxy,}{x+Dx^2+Fy^2.}
\end{equation*}
It is proved in~\cite{GGV04} that if $B=0$ then the period function of the center at the origin is globally monotone. So, from the point of view of the study of the period function, the most interesting stratum of quadratic centers is $B\neq 0,$ which can be brought by means of a rescaling to $B=1,$ i.e., 
\begin{equation}\label{Loud2}
 X_{\mu}\;\sist{-y+xy,}{x+Dx^2+Fy^2,}
\end{equation}
where $\mu\!:=(D,F)\in\R^2.$ This paper is addressed to study the period function of the center at the origin in this two-parametric family. More precisely, for a given $\hat\mu\in\R^2,$ we are concerned with the number of critical periodic orbits of $X_{\mu}$ with  $\mu\approx\hat\mu$ that emerge or disappear from the polycycle $\Pi_{\hat\mu}$ of $X_{\hat\mu}$ at the boundary of its period annulus as we move slightly the parameter. We refer to this number as the \emph{criticality} of the polycycle, $\mathrm{Crit}\left((\Pi_{\hat\mu},X_{\hat\mu}),X_{\mu}\right)$, see \defic{defi:criticality}. (Again this is the counterpart to the notion of cyclicity of a limit periodic set, see for instance \cite{Roussarie}.) Then we say that $\hat\mu\in\R^2$ is a local regular parameter if $\mathrm{Crit}\left((\Pi_{\hat\mu},X_{\hat\mu}),X_{\mu}\right)=0$ and that it is a local bifurcation parameter if $\mathrm{Crit}\left((\Pi_{\hat\mu},X_{\hat\mu}),X_{\mu}\right)\geqslant 1$. The initial work on this issue is \cite{MMV2} and, since our result is closely related to the ones obtained there, next we explain them succintly. With this aim, let~$\Gamma_U$ be the union of dotted straight
lines in \figc{diagrama}, whatever its colour is. Consider also the bold curve~$\Gamma_B.$ (Here the
subscripts~$B$ and~$U$ stand for bifurcation and unspecified respectively.)
 \begin{figure}[t]
 \centering
 \begin{lpic}[l(0mm),r(0mm),t(0mm),b(5mm)]{dib2(1)}
   \lbl[l]{7.5,73.5;$\mu_{\star}$}
   \lbl[l]{45.0,50.5;$D=\mathcal G(F)$}
  \end{lpic}     
\caption{Bifurcation diagram of the period function at the polycycle according to~\cite{MMV2} and, in colour, the subsequent improvements due to \cite{TopaV,MMSV,MV,Jordi}, where $\mu_{\star}=(-F_{\star},F_{\star})$ with $F_{\star}\approx 2.34.$ The curve that joins $\left(-\frac{3}{2},\frac{3}{2}\right)$ and $\left(-\frac{1}{2},1\right)$ is the graphic of an analytic function $D=\mathcal G(F).$}\label{diagrama}
 \end{figure}
Then, following this notation, \cite[Theorem A]{MMV2} shows that the open set $\R^2\setminus\left(\Gamma_B\cup\Gamma_U\right)$ corresponds to local regular parameters and that the ones in $\Gamma_B$ are local bifurcation parameters. The authors also conjecture that any parameter in $\Gamma_U$ is regular, except for the segment $\{0\}\!\times\!\left[0,\frac{1}{2}\right]$ in the vertical axis, that should consist of bifurcation parameters. The proof of the result in \cite{MMV2} is based 
on the explicit computation of the first non-vanishing term of the asymptotic expansion of the derivative of the period function. One of the major difficulties to tackle~$\Gamma_U$ lies in the necessity that this expansion is uniform with respect to $\mu$ in order to show that some parameter is regular. There have been however some progress in the study of the parameters in $\Gamma_U$: 
\begin{itemize}
\item From the results in \cite{TopaV,Jordi} it follows that the parameters in blue are indeed regular. In these papers the 
         authors determine a region $M$ in the parameter plane for which the corresponding center has a globally
          monotonic period function. The parameters that we draw in blue are inside the interior of~$M$, which prevents 
          the bifurcation of critical periodic orbits.
\item By \cite[Theorem C]{MMSV} it follows that the parameters in green are regular as well. In that paper the 
         authors give an asymptotic expansion of the Dulac time (time of the Dulac map) of an unfolding of a saddle-node.
         The techniques used in \cite{MMV2} enable only to study an unfolding of a hyperbolic saddle.  
\item Theorem B in \cite{MV} shows that the parameters in red, more precisely the segment 
          $\{0\}\!\times\!\left[\frac{1}{4},\frac{1}{2}\right]$, are bifurcation values of the period function at the polycycle. 
          To this end, after blowing up the 
         polycycle, the authors show that any neighbourhood of a parameter $\hat \mu\in \{0\}\!\times\!\left[\frac{1}{4},
         \frac{1}{2}\right]$ contains two parameters, say $\mu_+$ and $\mu_-$, such that the derivative of the 
         period function near the polycycle is positive for $X_{\mu_+}$ and negative for $X_{\mu_-}.$ 
\end{itemize} 
 \begin{figure}[t]
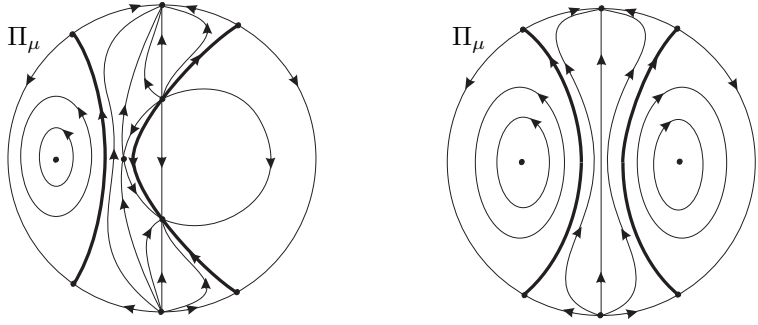

 \centering
 \begin{lpic}[l(0mm),r(0mm),t(0mm),b(5mm)]{dib3(2)}
   \lbl[l]{0,18.5;$\Pi_{\mu}$}
   \lbl[l]{29.25,18.5;$\Pi_{\mu}$}
  \end{lpic}     
\caption{\label{disc}
Phase portrait of \refc{Loud2} in the PoincarŽ disc for $\mu=(D,F)\in\Lambda$ with $D<-1$ (left) and $D>-1$ (right), where for convenience we place the center at $(0,0)$ on the left of the centered invariant line $\{x=1\}.$ The polycycle $\Pi_{\mu}$ at the outer boundary of the period annulus is the same in both cases: two hyperbolic saddles at infinity and the heteroclinic orbits between them. The invariant hyperbola $\mathscr C_{\mu}$ is in boldface type.}
\end{figure}
Beyond the dichotomy regular vs bifurcation, a challenging problem not tackled so far is the study of the criticality of the local bifurcation parameters, i.e., to compute the \emph{exact number} of critical periodic orbits that bifurcate from the polycycle. The present paper is addressed to this problem. 

Besides studying the local problem explained above, the authors in~\cite{MMV2} obtained a number of results concerning the global behaviour of the period function of~\eqref{Loud2}. The combination of all these results with the ones obtained by Chicone and Jacobs \cite{ChiJac}, lead them to formulate a very precise conjecture for the bifurcation diagram of the (global) period function. According to this conjecture, the criticality of the polycycle should be equal to one for any $\mu\in\Gamma_{B}\setminus\{(0,0),(-\frac{3}{2},\frac{3}{2}),(-2,2)\},$ whereas it should be equal to two for the three remaining parameters. Note that~$\Gamma_B$ is the union of some explicit straight segments and a curve that joins the points $\left(-\frac{3}{2},\frac{3}{2}\right)$ and $\left(-\frac{1}{2},1\right),$ which can be proved to be the graphic of an analytic function $D=\mathcal G(F),$ see \propc{prop:corba_bif}. Our main result shows that the criticality in a portion of this curve is exactly one. For the sake of completeness we also reprove some results already obtained in \cite{MMV2}. The next statement gathers these results and we stress that it concerns the parameters inside
 \begin{equation}\label{Lambda}
 \Lambda\!:=\{(D,F)\in\R^2: F>1, D<-1/2, D+F>0\}.
 \end{equation}
Accordingly, see \figc{diagrama}, the above mentioned curve is inside $\Lambda$. On the other hand, if $(D,F)\in\Lambda\cap\Gamma_B$ then $(F-2)(D-\mathcal{G}(F))=0.$ 

\begin{thmx}\label{thm:crit}
Let us consider the period function of the center at the origin of the differential system~\eqref{Loud2} for $\mu\in\Lambda$. Then the open set $\Lambda\setminus\left(\Gamma_B\cup\{F=3/2\}\right)$ corresponds to local regular values of the period function at the outer boundary of the period annulus. On the other hand, the parameters in $\Gamma_B$ are local bifurcation values of the period function at the outer boundary of the period annulus. Moreover $\mathrm{Crit}\left((\Pi_{\hat\mu},X_{\hat\mu}),X_{\mu}\right)=1$ for any $\hat\mu=(\hat D, \hat F)$ with $\hat D=\mathcal{G}(\hat F)$ and $\frac{4}{3}<\hat F<\frac{3}{2}$. 
\end{thmx}

\figc{disc} displays the phase portrait of \refc{Loud2} for the parameter values under consideration in \teoc{thm:crit}. We remark that the new result, and the main goal of this paper, is the last assertion. To prove \teoc{thm:crit} we apply the tools developed in \cite{ManRojVil2016-JDE,ManRojVil2016-JDDE}, which provide sufficient conditions in order that $\mathrm{Crit}\left((\Pi_{\hat\mu},X_{\hat\mu}),X_{\mu}\right)\leqslant n.$ The underlying idea of these conditions is to guarantee that the derivative of the period function can be embedded in an extended complete Chebyshev system of dimension $n+1$ in a neighbourhood of the polycycle. This is a completely different approach from the previous works \cite{MMV2,MMSV,MV}, which rely on the use of normal forms near the singularities of the polycycle. At this respect let us point out that we recover all the results in \cite{MMV2} regarding the dichotomy regular vs bifurcation in the subset~$\Lambda$ except for the regularity of the segment in $\{F=\frac{3}{2}\}.$ At these parameters, the asymptotic expansion of the period function has a logarithmic term, which forces the use of the so-called \emph{Roussarie-Ecalle compensator} in order to guarantee the necessary uniformity with respect to $\mu$ (see \cite[Theorem 3.6]{MMV2}). The tools we have developed so far do not allow to deal with this scenario. Let us also remark that the tools in \cite{ManRojVil2016-JDE,ManRojVil2016-JDDE} are in fact addressed to potential differential systems, which is not the case of~\eqref{Loud2}. We avoid this problem by appealing to \cite[Lemma~14]{GGV04}, that gives a class of integrable differential systems that can be brought to a potential system by means of an explicit coordinate transformation. Luckily the differential system under consideration~\eqref{Loud2} is inside this class.   

For the sake of completeness, let us finish this section by quoting some other results about the period function of~\eqref{Loud2}. Thus, by applying mainly techniques based on Picard-Fuchs equations for algebraic curves, Yulin Zhao {\it et al} proved in~\cite{Zhao2,Zhao3,Zhao5,Zhao4} that system~\eqref{Loud2} has at most one critical periodic orbit for any $\mu=(D,F)$ with $F\in\left\{\frac{1}{4},\frac{3}{4},\frac{3}{2},2\right\}$. On the other hand, R. Chouikha showed in \cite{Chouikha} the monotonicity of the period function for the parameters in the straight lines $F+2D=1$ and $F=-1$ and some other segments inside $D=\frac{1}{2}$, $D=0$, $F=1$ and $F= 2$.

The paper is organized in the following way. In \secc{defi} we introduce the related definitions and we explain the previous results obtained in \cite{ManRojVil2016-JDE,ManRojVil2016-JDDE} that we shall use in the proof of \teoc{thm:crit}, which is carried out in \secc{provaA}.  

\section{Definitions and previous results}\label{defi}

A singular point $p$ of an analytic vector field $X=f(x,y)\partial_x+g(x,y)\partial_y$ is a \emph{center} if it has a punctured neighbourhood that consist entirely of periodic orbits surrounding $p$. The largest punctured neighbourhood with this property is called the \emph{period annulus} of the center and henceforth it will be denoted by $\PA$. From now on $\partial\PA$ will denote the boundary of $\PA$ after embedding it into $\RP^2$. Clearly the center $p$ belongs to $\partial\PA$, and in what follows we will call it the \emph{inner boundary} of the period annulus. We also define the \emph{outer boundary} of the period annulus to be $\Pi\!:=\partial\PA\setminus\{p\}$. Note that $\Pi$ is a non-empty compact subset of~$\RP^2$. In case that~$X$ is polynomial then, by means of the PoincarŽ compactification,  it has a meromorphic extension $\hat X$ to infinity and the outer boundary turns out to be a polycycle of $\hat X.$ The \emph{period function} of the center assigns to each periodic orbit in $\PA$ its period. Since the period function is defined on the set of periodic orbits in~$\mathscr P,$ in order to study its qualitative properties usually the first step is to parametrize this set. This can be done by taking an analytic transverse section to $X$ on~$\mathscr P$, for instance an orbit of the orthogonal vector field~$X^{\perp}$. If $\{\gamma_s\}_{s\in (a,b)}$ is such a parametrization, then $s\longmapsto T(s)\!:=\!\{\mbox{period of $\gamma_s$}\}$ is an analytic map that provides the qualitative properties of the period function that we are concerned about. In particular the existence of \emph{critical periods}, which are isolated critical points of this function, i.e. $\hat s\in (a,b)$ such that $T'(s)=\alpha(s-\hat s)^k+\op\bigl((s-\hat s)^k\bigr)$ with $\alpha\neq 0$ and $k\geqslant 1.$ In this case we shall say that $\gamma_{\hat s}$ is a \emph{critical periodic orbit} of multiplicity $k$ of the center. One can readily see that this definition does not depend on the particular parametrization of the set of periodic orbits used.


Our aim in this paper is to study the bifurcation of critical periodic orbits from the outer boundary of the period annulus. As any bifurcation phenomenon, this occurs in case that $X$ depends on a parameter, say $\mu\in\Lambda\subset\R^d$. Thus, for each $\mu\in\Lambda$, assume that $X_{\mu}$ is an analytic vector field on some open set $\mathscr U_{\mu}$ of~$\R^2$ with a center at $p_{\mu}$. Following the notation introduced before, we denote by $\Pi_{\mu}$ the outer boundary of its period annulus $\PA_{\mu}$. Concerning the regularity with respect to $\mu$, we assume that $\{X_{\mu}\}_{\mu\in\Lambda}$ is a continuous family, i.e., that the map $(x,y,\mu)\longmapsto X_{\mu}(x,y)$ is continuous in $\{(x,y,\mu)\in\R^{d+2}: (x,y)\in\mathscr U_{\mu}, \mu\in\Lambda\}$.
In the following definition $d_H$ stands for the Hausdorff distance between compact sets of $\RP^2$.

\begin{defi}\label{defi:criticality}
Consider a continuous family $\{X_{\mu}\}_{\mu\in\Lambda}$ of planar analytic vector fields with a center and fix some $\hat\mu\in\Lambda$. Suppose that the outer boundary of the period annulus varies continuously at $\hat\mu\in\Lambda$, meaning that for any $\epsilon>0$ there exists $\delta>0$ such that $d_H(\Pi_{\mu},\Pi_{\hat\mu})\leq\epsilon$ for all $\mu\in\Lambda$ with $\|\mu-\hat\mu\|\leq\delta$. Then, setting
\[
N(\delta,\epsilon)=\sup\{\# \text{ critical periodic orbits }\gamma\text{ of }X_{\mu}\text{ in }\PA_{\mu}\text{ with }d_H(\gamma,\Pi_{\hat\mu})\leq\epsilon\text{ and }\|\mu-\hat\mu\|\leq\delta\},
\]
the \emph{criticality} of $(\Pi_{\hat\mu},X_{\hat\mu})$ with respect to the deformation $X_{\mu}$ is $\text{Crit}\bigl((\Pi_{\hat\mu},X_{\hat\mu}),X_{\mu}\bigr)\!:=\inf_{\delta,\epsilon}N(\delta,\epsilon)$.
\end{defi}
In the previous definition $\text{Crit}\bigl((\Pi_{\hat\mu},X_{\hat\mu}),X_{\mu}\bigr)$ may be infinite but if it is not, then it gives the maximal number of critical periodic orbits of $X_{\mu}$ that tend to $\Pi_{\hat\mu}$ in the Hausdorff sense as $\mu\rightarrow\hat\mu$.

\begin{defi}
We say that $\hat\mu\in\Lambda$ is a local regular value of the period function at the outer boundary of the period annulus if $\text{Crit}\bigl((\Pi_{\hat\mu},X_{\hat\mu}),X_{\mu}\bigr)=0$. Otherwise we say that it is a local bifurcation value of the period function at the outer boundary.
\end{defi}

We shall next state the results from \cite{ManRojVil2016-JDE,ManRojVil2016-JDDE} that we shall use in the proof of \teoc{thm:crit}. These papers are concerned with analytic potential differential systems
\[
Y_{\mu}\;\begin{cases}\dot{x}=-y,\\ \dot{y}=V_{\mu}'(x),\end{cases}
\]
depending on a parameter $\mu\in\Lambda\subset\R^d.$ Here, for each fixed $\mu\in\Lambda,$ $V_{\mu}$ is an analytic function on a certain real interval $I_{\mu}$ that contains $x=0$. We shall also use the vector
field notation $Y_{\mu}=-y\partial_x+V'_{\mu}(x)\partial_y$ to
refer to the above differential system. We suppose $V'_{\mu}(0)=0$ and $V''_{\mu}(0)>0$, so that the origin is a 
non-degenerated center and we shall denote the
projection of its period annulus~$\PA_{\mu}$ on the $x$-axis by
$\PI_{\mu}=(x_{\ell}(\mu),x_r(\mu))$. Thus
$x_{\ell}(\mu)<0<x_r(\mu).$ The corresponding Hamiltonian function
is given by $H_{\mu}(x,y)=\frac{1}{2}y^2+V_{\mu}(x),$ where we fix that
$V_{\mu}(0)=0,$ and we set the
energy level of the outer boundary of $\PA_{\mu}$ to be $h_0(\mu)$,
i.e. $H_{\mu}(\PA_{\mu})=(0,h_0(\mu)).$ Note then that $h_0(\mu)$ is a
positive number or $+\infty.$ In addition we define
\[
 g_{\mu}(x)\!:=x\sqrt{\frac{V_{\mu}(x)}{x^2}}=\mathrm{sgn}(x)\sqrt{V_{\mu}(x)},
\]
which is clearly a diffeomorphism from $(x_{\ell}(\mu),x_r(\mu))$ to $(-\sqrt{h_0(\mu)},\sqrt{h_0(\mu)})$ due to $V_{\mu}(0)=V'_{\mu}(0)=0$ and $V''_{\mu}(0)>0.$ In order to state the above mentioned results appropiately, it is necessary to introduce a number of definitions:

\begin{defi}\label{def:Ak}
We say that the family of potential analytic differential systems $\{Y_{\mu}\}_{\mu\in\Lambda}$ verifies the hypothesis \hypertarget{H}{\hyperlink{H}{\textnormal{\textbf{(H)}}}} in case that:
\begin{enumerate}[$(a)$]
\item For all $k\geq 0$, the map $(x,\mu)\longmapsto V_{\mu}^{(k)}(x)$ is continuous on $\{(x,\mu)\in\R\times\Lambda : x\in I_{\mu}\},$
\item $\mu\longmapsto x_{r}(\mu)$ is continuous on $\Lambda$ or $x_{r}(\mu)=+\infty$ for all $\mu\in\Lambda,$
\item $\mu\longmapsto x_{\ell}(\mu)$ is continuous on $\Lambda$ or $x_{\ell}(\mu)=-\infty$ for all $\mu\in\Lambda,$
\item $\mu\longmapsto h_0(\mu)$ is continuous on $\Lambda$ or $h_0(\mu)=+\infty$ for all $\mu\in\Lambda.$
\end{enumerate}\vspace*{-.63cm}
\end{defi}

\begin{defi}
Let $f$ be an analytic function on $(a,b)$. We say that $f$ is \emph{quantifiable} at $b$ by $\alpha$ with limit~$\ell$ in case that:
\begin{enumerate}[$(a)$]
\item If $b\in\R$, then $\lim_{x\rightarrow b^-}f(x)(b-x)^{\alpha}=\ell$ and $\ell\neq 0$.
\item If $b=+\infty$, then $\lim_{x\rightarrow+\infty}f(x)x^{-\alpha}=\ell$ and $\ell\neq 0$.
\end{enumerate}
We call $\alpha$ the \emph{quantifier} of $f$ at $b$. We shall use the analogous definition at $a$.
\end{defi}


\begin{defi}
Let $\{f_{\mu}\}_{\mu\in\Lambda}$ be a continuous family of analytic functions on $(a(\mu),b(\mu))$. Assume that $b$ is either a continuous function from $\Lambda$ to $\R$ or $b(\mu)=+\infty$ for all $\mu\in\Lambda$. Given $\hat\mu\in\Lambda$ we shall say that $\{f_{\mu}\}_{\mu\in\Lambda}$ is \emph{continuously quantifiable} in $\hat\mu$ at $b(\mu)$ by $\alpha(\mu)$ with limit $\ell(\hat\mu)$ if there exists an open neighbourhood $U$ of~$\hat\mu$ such that $f_{\mu}$ is quantifiable at $b(\mu)$ by $\alpha(\mu)$ with limit $\ell(\mu)$ for all $\mu\in\Lambda$ and, moreover,
\begin{enumerate}[$(a)$]
\item If $b(\hat\mu)<+\infty$, then $\lim_{(x,\mu)\rightarrow(b(\hat\mu),\hat\mu)}f_{\mu}(x)(b(\mu)-x)^{\alpha(\mu)}=\ell(\hat\mu)$ and $\ell(\hat\mu)\neq 0$.
\item If $b(\hat\mu)=+\infty$, then $\lim_{(x,\mu)\rightarrow(+\infty,\hat\mu)}f_{\mu}(x)x^{-\alpha(\mu)}=\ell(\hat\mu)$ and $\ell(\hat\mu)\neq 0$.
\end{enumerate}
For the sake of shortness, in the first case we shall write $f_{\mu}(x)\sim_{b(\mu)} \ell(\mu)(b(\mu)-x)^{-\alpha(\mu)}$ at $\hat\mu$, and in the second case $f_{\mu}(x)\sim_{+\infty} \ell(\mu)x^{\alpha(\mu)}$ at $\hat\mu$.
\end{defi}

\begin{defi}\label{wronskian}
Let $f_0,f_1,\dots,f_{k-1}$ be analytic functions on an open interval $I$ of $\R$. Then
\[
W[f_0,f_1,\dots,f_{k-1}](x)=\text{det}\left(f_j^{(i)}(x)\right)_{0\leq i,j\leq k-1}=
\left|\begin{array}{ccc}
f_0(x) & \cdots & f_{k-1}(x)\\
f_0'(x) & \cdots & f_{k-1}'(x)\\
		&	\vdots & \\
f_0^{(k-1)}(x) & \cdots & f_{k-1}^{(k-1)}(x)
\end{array}\right|
\]
is the \emph{Wronskian} of $(f_0,f_1,\dots,f_{k-1})$ at $x\in I$.
\end{defi}

\begin{defi}
Given $\nu_1,\dots,\nu_n\in\R$, we consider the linear ordinary differential operator 
\[
\DD_{\boldsymbol\nu_n}\!:\mathscr C^{\omega}\bigl((0,1)\bigr)\longrightarrow \mathscr C^{\omega}\bigl((0,1)\bigr)
\]
defined by
\begin{equation*}
\DD_{\boldsymbol\nu_n}[f](x)\!:=(x(1-x^2))^{\frac{n(n+1)}{2}}\frac{W\left[\psi_{\nu_1},\ldots,\psi_{\nu_n},f\right](x)}{\prod_{i=1}^n\psi_{\nu_i}(x)},
\end{equation*}
where $\psi_{\nu}(x)\!:=\frac{1}{1-x^2}\left(\frac{x}{\sqrt{1-x^2}}\right)^{\nu}$ and, for shortness, we use the notation $\boldsymbol\nu_n=(\nu_1,\dots,\nu_n)$. In addition we define $\DD_{\boldsymbol\nu_0}\!:=id$ for the sake of convenience.
\end{defi}

\begin{defi}\label{def_moment}
Let $f$ be an analytic function on $[0,1).$ Then, for each $n\in\N$ we call
\[
N_n[f]:=\int_0^{1} \frac{f(x)}{\sqrt{1-x^2}}\left(\frac{x}{\sqrt{1-x^2}}\right)^{2n-2}dx
\]
the $n$-th \emph{momentum} of $f$, whenever it is well defined.
\end{defi}

We are now in position to state \cite[Theorem B]{ManRojVil2016-JDDE}, which constitutes the main ingredient in the proof of \teoc{thm:crit}. In its statement, the assumptions requiring the existence of functions $\nu_1,\nu_2,\ldots,\nu_n$ and that $N_1\equiv N_2\equiv \ldots\equiv N_{j-1}\equiv 0$ must be considered void for $n=0$ and $j=1$, respectively. Moreover, for a given function $f$, we use the notation $\mathcal P[f](x)\!:=f(x)+f(-x)$.

\begin{thm}\label{teorema-jdde}
Let $\{Y_{\mu}\}_{\mu\in\Lambda}$ be a family of potential analytic systems verifying \hyperlink{H}{\textnormal{\textbf{(H)}}} such that $h_0(\mu)<+\infty$ for all $\mu\in\Lambda$.  Assume that there exist $n\geq 0$ continuous functions $\nu_1,\nu_2,\ldots,\nu_n$ in a neighbourhood of some fixed $\hat\mu\in\Lambda$ such that the family 
\begin{equation}\label{familia}
\bigl\{(\DD_{\boldsymbol\nu_n(\mu)}\circ\mathcal P)\bigl[z\sqrt{h_0(\mu)}(g_{\mu}^{-1})''(z\sqrt{h_0(\mu)})\bigr]\bigr\}_{\mu\in\Lambda}
\end{equation}
is continuously quantifiable in $\Lambda$ at $z=1$ by
$\xi(\mu).$ For each $i\in\N$, let $N_i(\mu)$ be the $i$-th momentum of $(\DD_{\boldsymbol\nu_n(\mu)}\circ\mathcal P)\bigl[z\sqrt{h_0(\mu)}(g_{\mu}^{-1})''(z\sqrt{h_0(\mu)})\bigr]$, whenever it is well defined. The following assertions hold:
\begin{enumerate}[$(a)$]
\item If $\xi(\hat\mu)>\frac{1}{2}$, then $\mathrm{Crit}\bigl((\out_{\hat\mu},Y_{\hat\mu}),X_{\mu}\bigr)\leqslant n$.
\item If $\xi(\hat\mu)<\frac{1}{2}$, let $m\in\N$ be such that $\xi(\hat\mu)+m\in \left[\frac{1}{2},\frac{3}{2}\right)$. Then $\mathrm{Crit}\bigl((\out_{\hat\mu},Y_{\hat\mu}),X_{\mu}\bigr)\leqslant n$ in case that
\begin{enumerate}[$(b1)$]
\item either $N_1\equiv N_2\equiv \ldots\equiv N_{j-1}\equiv 0$ and $N_j(\hat\mu)\neq 0$ for some $j\in\{1,2,\ldots,m\}$,
\item or $N_1\equiv N_2\equiv \ldots\equiv N_m\equiv 0$ and  $\xi(\hat\mu)+m\notin\left\{\frac{1}{2},1\right\}.$
\end{enumerate}
\end{enumerate}
Finally, if the following conditions are verified, then the family \refc{familia}  is continuously quantifiable at $z=1$ by
$\xi(\mu)=-\min\bigl\{\bigl(\frac{\alpha_{\ell}}{\beta_{\ell}}\bigr)(\mu), \bigl(\frac{\alpha_r}{\beta_r}\bigr)(\mu)\bigr\}-\frac{1}{2}\sum_{i=1}^n\nu_i(\mu)-\frac{n(n+1)}{2}+1:$
\begin{enumerate}[$(i)$]
\item $\{h_0(\mu)-V_{\mu}\}_{\mu\in\Lambda}$ is continuously quantifiable at $x_{\ell}(\mu)$ by
$\beta_{\ell}(\mu)$ and at $x_{r}(\mu)$ by $\beta_r(\mu)$ with limits $b_{\ell}(\mu)$ and $b_r(\mu)$, respectively,
\item setting $\mathscr R_{\mu}\!:=\frac{(V_{\mu}')^2-2V_{\mu}V_{\mu}''}{(V_{\mu}')^3},$ the function 
\[
x\longmapsto 
  V_{\mu}'(x)^{-\frac{n(n+1)}{2}}W\!\!\left[\left(\frac{V_{\mu}}{h_0(\mu)-V_{\mu}}\right)^{\frac{1}{2}\nu_1(\mu)},\ldots,\left(\frac{V_{\mu}}{h_0(\mu)-V_{\mu}}\right)^{\frac{1}{2}\nu_n(\mu)},(h_0(\mu)-V_{\mu})\,{V_{\mu}}^{\!\frac{1}{2}}\,\mathscr R_{\mu}\right]\!\!(x)
\]
is continuously quantifiable at $x_{\ell}(\mu)$ by $\alpha_{\ell}(\mu)$ and at $x_{r}(\mu)$
by $\alpha_r(\mu)$ with limits $a_{\ell}(\mu)$ and $a_r(\mu)$,
respectively,
\item and either $\frac{\alpha_{\ell}}{\beta_{\ell}}(\mu)\neq\frac{\alpha_r}{\beta_r}(\mu)$  or,
otherwise,
$\bigr(a_{r}(b_r)^{-\frac{\alpha_r}{\beta_r}}+(-1)^{\frac{n(n+1)}{2}}a_{\ell}(b_{\ell})^{-\frac{\alpha_{\ell}}{\beta_{\ell}}}\bigl)(\mu)\neq 0.$
\end{enumerate}
\end{thm}

As we already mentioned, the idea behind this result is to give conditions in order that the derivative of the period function can be embedded in an extended complete Chebyshev system of dimension $n+1$ in a neighbourhood of the polycycle. More specifically, denoting the period of the periodic orbit of $Y_{\mu}$ inside the energy level $\{H_{\mu}(x,y)=h\}$ by $T_{\mu}(h),$
these conditions guarantee that
\[
\lim_{z\longrightarrow 1}(1-z)^{\nu_n(\mu)}W\bigl[\psi_{\nu_1(\mu)}(z),\dots,\psi_{\nu_{n-1}(\mu)}(z),T_{\mu}'(z^2h_0(\mu))\bigr]=\Delta^{\!\star}(\mu),
\]
uniformly in $\mu\approx \hat\mu$, and that $\Delta^{\!\star}(\hat\mu)\neq 0.$ At this respect, the following observation will enable us to avoid some cumbersome computations. 

\begin{rem}\label{remark}
From the proof of \teoc{teorema-jdde} in \cite{ManRojVil2016-JDDE} for the particular case $n=0$ it follows that 
 \[
  \frac{T_{\mu}'(h)}{(h_0(\mu)-h)^{\gamma(\mu)}}\longrightarrow \Delta^{\!\star}(\mu)\text{ as $(h,\mu)\longrightarrow (h_0(\hat\mu),\hat\mu)$,}
 \]
where
\[
\begin{cases}
\gamma(\mu)= \frac{1}{2}-\xi(\mu)  \text{ and } \Delta^{\!\star}(\mu)= C(\mu)\delta(\mu) & \text{in cases }(a)\text{ and }(b2),\\[5pt]
\gamma(\mu)= 1-j \text{ and } \Delta^{\!\star}(\mu)= C(\mu)N_j(\mu) & \text{in case }(b1),
\end{cases}
\]
with $C(\mu)>0$ for all $\mu\approx\hat\mu$, and
\[
\delta(\mu)=\begin{cases}
a_r(b_r)^{-\frac{\alpha_r}{\beta_r}} & \text{ if }\frac{\alpha_r}{\beta_r}<\frac{\alpha_{\ell}}{\beta_{\ell}},\\
a_r(b_r)^{-\frac{\alpha_r}{\beta_r}}+a_{\ell}(b_{\ell})^{-\frac{\alpha_{\ell}}{\beta_{\ell}}} & \text{ if }\frac{\alpha_r}{\beta_r}=\frac{\alpha_{\ell}}{\beta_{\ell}},\\
a_{\ell}(b_{\ell})^{-\frac{\alpha_{\ell}}{\beta_{\ell}}} & \text{ if }\frac{\alpha_r}{\beta_r}>\frac{\alpha_{\ell}}{\beta_{\ell}}.\\
\end{cases}
\]
%
\end{rem}

We shall also apply the following technical result (see \cite[Lemma 3.12]{ManRojVil2016-JDDE}).

\begin{lema}\label{lem:reduc_moment_finit}
Let $f$ be an analytic function on $[0,1)$, $\nu_1,\nu_2,\ldots,\nu_{n}\in\R$ and $\ell\in\N.$  Let us assume that $\DD_{\boldsymbol\nu_{n-1}}[f]$ is quantifiable at $1$ by $\xi$. If $\xi<3/2-\ell$, then
\[
N_{\ell}\bigl[\DD_{\boldsymbol\nu_{n}}[f]\bigr]=c_n(1-2\ell-\nu_n)N_{\ell}\bigl[\DD_{\boldsymbol\nu_{n-1}}[f]\bigr],
\]
where $c_1\!:=1$ and $c_n\!:=\prod_{i=1}^{n-1}(\nu_{n}-\nu_i)$ for $n\geqslant 2.$
\end{lema}

\section{Proof of Theorem~A}\label{provaA}

We point out that in this section $\Lambda$ refers to the parameter subset as defined in \refc{Lambda}. It is well known, see for instance \cite{MMV2}, that if $F\notin\{0,1,\frac{1}{2}\}$
then the differential system \refc{Loud2} has a first integral given by
\begin{equation}\label{first-integral}
\mbox{$H_{\mu}(x,y)=(1-x)^{-2F}\left(\frac{1}{2}\,y^2-q_{\mu}(x)\right)$,}
\end{equation}
where $q_{\mu}(x)=a(\mu)x^2+b(\mu)x+c(\mu)$ with
\begin{equation*}
a=\frac{D}{2(1-F)},\mbox{ }
b=\frac{D-F+1}{(1-F)(1-2F)}\mbox{ and }
c=\frac{F-D-1}{2F(1-F)(1-2F)}.
\end{equation*}
Its corresponding integrating factor is $\kappa(x)=(1-x)^{-2F-1}$. 
In addition, the line at infinity $L_{\infty}$, the conic
$\mathscr{C}_{\mu}=\{\frac{1}{2}y^2-q_{\mu}(x)=0\}$ and the line $\{x=1\}$
are invariant curves of \refc{Loud2}. If $\mu\in\Lambda$ then $\mathscr C_{\mu}$ is a hyperbola that intersects $y=0$ at 
 \[
x=p_{1}(\mu)\!:=\frac{-b-\sqrt{b^2-4ac}}{2a}\text{ and }x=p_{2}(\mu)\!:=\frac{-b+\sqrt{b^2-4ac}}{2a},
\]
with $0<p_1(\mu)<p_2(\mu)$. Moreover, see \figc{disc}, for these parameter values, the outer boundary of the period annulus of the center at the origin consists of the branch of the hyperbola $\mathscr C_{\mu}$ passing through the point $(p_1,0)$ and the line at infinity $L_{\infty}$ joining two hyperbolic saddles. The next result gathers some relevant facts proved in \cite{MMV2} that we shall use later on. 

\begin{prop}\label{prop:corba_bif}
For each $\mu\in\Lambda$ and $s\approx 0$ positive, let $P(s;\mu)$ be the period of the periodic orbit of $\refc{Loud2}$ passing through the point $\bigl(p_1(\mu)-s,0\bigr).$ If $F\in (1,\frac{3}{2})$ then $\lim_{s\to 0^+}P'(s;\mu)=\Delta(\mu)$ uniformly on compact subsets of $\Lambda$, where
 \[
  \Delta(\mu)=\frac{-1/\sqrt{2a}}{(p_2-p_1)(1-p_1)}\left\{2-\int_0^1\left(u^{2(1-F)}\left(\frac{1-p_2}{1-p_1}\,(u-1)+1\right)^{2F-1}-1
             \right)\frac{du}{(1-u)^{3/2}}\right\}.
 \]
Moreover, the set $\{\mu\in\Lambda:\Delta(\mu)=0\}$ is the graph of an analytic function $D=\mathcal G(F)$ defined for $F\in \left(1,\frac{3}{2}\right)$ and satisfying the following properties:
\begin{enumerate}[$(a)$]
 \item $-F<\mathcal G(F)<-\frac{1}{2}$ for all $F\in \left(1,\frac{3}{2}\right)$, 
 \item $\lim_{F\to\frac{3}{2}} \mathcal G(F)=-\frac{3}{2}$, and
 \item $\lim_{F\to 1} \mathcal G(F)=-\frac{1}{2}.$
\end{enumerate}
\end{prop}

\begin{prova}
The first assertion follows from $(a)$ in \cite[Theorem 3.6]{MMV2} and the properties of the set $\{\Delta(\mu)=0\}$ from \cite[Proposition 3.11]{MMV2}.
\end{prova}

Note in particular that $\{\Delta(\mu)=0\}$ is an analytic curve inside $\Lambda$ that joins the parameters $\mu=\left(-\frac{3}{2},\frac{3}{2}\right)$ and $\mu=\left(-\frac{1}{2},1\right).$ Recall at this point that the differential system \refc{Loud2} has a first integral which is quadratic in~$y$, see \refc{first-integral}, and that its corresponding integrating factor depends only on $x.$ Taking advantage of this, by applying \cite[Lemma~14]{GGV04} it follows that the coordinate transformation 
\begin{equation*}
(u,v)=\bigl(\phi(1-x), (1-x)^{-F}y\bigr),\text{ where $\phi(z)\!:=\frac{z^{-F}-1}{F}$},
\end{equation*} 
brings the differential system \eqref{Loud2} to the potential system
\begin{equation}\label{Loud-potential}
\begin{cases}
\dot{u} = -v,\\[2pt]
\dot{v} = (Fu+1)\bigl((Fu+1)^{-\frac{1}{F}}-1\bigr)\bigl(D(Fu+1)^{-\frac{1}{F}}-D-1\bigr).
\end{cases}
\end{equation}
Evidently this differential system has a non-degenerated center at the origin. The projection on the $u$-axis of its period annulus 
is the interval $\mathcal I_{\mu}\!:=\left(-\frac{1}{F},u_r(\mu)\right)$, where
 \[
 u_r(\mu)\!:=\frac{(1-p_1(\mu))^{-F}-1}{F}.
\]
This interval is precisely the image by $x\longmapsto \phi(1-x)$ of $\bigl(-\infty,p_1\bigr)$. Let $H(u,v)=\frac{1}{2}v^2+V_{\mu}(u)$ be the Hamiltonian function of the potential system \refc{Loud-potential} with $V_{\mu}(0)=0.$ Setting $z=\phi^{-1}(u)=(Fu+1)^{-1/F}$ for shortness, one can check that 
 \begin{equation}\label{eq:expressions_V}
 \begin{array}{ll}
  V_{\mu}(u)=h_0(\mu)-z^{-2F}V_0(z,\mu) & \text{ with }V_0(z,\mu)=\frac{D}{2-2F}z^2+\frac{1+2D}{2F-1}z-\frac{D+1}{2F}, \\[3pt]
  V_{\mu}'(u)=z^{-F}V_1(z,\mu) & \text{ with }V_1(z,\mu)=(z-1)(D(z-1)-1),  \\[3pt]
  V_{\mu}''(u)=V_2(z,\mu) & \text{ with }V_2(z,\mu)=D(F-2)z^2-(2D+1)(F-1)z+F(D+1),  \\[3pt]
  V_{\mu}^{(3)}(u)=z^FV_3(z,\mu)& \text{ with }V_3(z,\mu)=-2D(F-2)z^2+(2D+1)(F-1)z,  \\[3pt]
   V_{\mu}^{(4)}(u)=z^{2F}V_4(z,\mu)& \text{ with }V_4(z,\mu)=2D(F^2-4)z^2-(2D+1)(F^2-1)z. 
 \end{array}
 \end{equation}
Here $h_0(\mu)\!:=\frac{F-D-1}{2F(F-1)(2F-1)}$ turns out to be the energy level of the outer boundary of the center at the origin for all $\mu\in\Lambda.$ 

In view of the properties explained above, the family of potential systems \refc{Loud-potential} with $\mu\in\Lambda$ satisfies the hypothesis \hyperlink{H}{\textnormal{\textbf{(H)}}} in \defic{def:Ak}. Moreover, since
this family is conjugated to \eqref{Loud2}, in order to prove \teoc{thm:crit} we can apply \teoc{teorema-jdde} to \refc{Loud-potential}. Our task now is to quantify all the functions involved in its application. This is the aim of the following lemmas. 

\begin{lema}\label{lema:traduccio}
Let us take $\mu\in\Lambda$. 
If $\bigl(f_{\mu}\circ\phi\bigr)(z)\sim_{+\infty} \ell(\mu)z^{\alpha(\mu)}$ 
then $f_{\mu}(u) \sim_{-\frac{1}{F}} \ell(\mu)(Fu+1)^{-\alpha(\mu)/F}$.
\end{lema}

\begin{prova}
Note that $\phi^{-1}(u)=(Fu+1)^{-1/F}\longrightarrow +\infty$ as $u$ tends to $-1/F$, uniformly on compact subsets of~$\Lambda.$ Let us fix $\hat\mu=(\hat D,\hat F)\in\Lambda.$ Then
\[
\lim_{(u,\mu)\to (-1/\hat F,\hat\mu)}f_{\mu}(u)\left(Fu+1\right)^{\frac{\alpha(\mu)}{F}}=
\lim_{(u,\mu)\to (-1/\hat F,\hat\mu)}\frac{\bigl(f_{\mu}\circ \phi\bigr)(\phi^{-1}(u))}{(\phi^{-1}(u))^{\alpha(\mu)}}=
\lim_{(z,\mu)\to (+\infty,\hat\mu)} \frac{\bigl(f_{\mu}\circ\phi\bigr)(z)}{z^{\alpha(\mu)}}=\ell(\hat\mu)
\]
and this shows the result. 
\end{prova}

\begin{lema}\label{lema:quantificar_V}
For all $\mu\in\Lambda$ the following holds:
\begin{enumerate}[$(a)$]
\item $h_0(\mu)-V_{\mu}(u) \sim_{-\frac{1}{F}} \frac{D}{2-2F} (Fu+1)^{2-2/F}$,
\item $V_{\mu}'(u_r(\mu))\neq 0$ and $h_0(\mu)-V_{\mu}(u) \sim_{u_r(\mu)} V_{\mu}'(u_r(\mu))(u_r(\mu)-u)$.
\end{enumerate} 
\end{lema}

\begin{prova}
From \refc{eq:expressions_V}, $h_0(\mu)-\bigl(V_{\mu}\circ\phi\bigr)(z)=z^{-2F}V_0(z,\mu)\sim_{+\infty} \frac{D}{2-2F}z^{2-2F}$. Then by applying~\lemc{lema:traduccio} the assertion in~$(a)$ follows. The analyticity of $V_{\mu}$ at $u=u_r$, together with the fact that $V_{\mu}'(u_r)\neq 0$ due to $V_1(1-p_1(\mu),\mu)=p_1(\mu)(Dp_1(\mu)+1)\neq 0$, easily imply~$(b)$. So the result is proved. 
\end{prova}

\begin{lema}\label{lema:corba}
If $\mu\in\Lambda$ then $V_{\mu}'(u)^2-2V_{\mu}(u)V_{\mu}''(u)$ is strictly positive
at $u=u_r(\mu).$ 
\end{lema}

\begin{prova}
Since $V(u_r)=h_0$, we must prove that $V_{\mu}'(u_r)^2-2h_0(\mu)V_{\mu}''(u_r)>0$ for all $\mu\in\Lambda.$ To this end, on  account of \refc{eq:expressions_V}, we write $V_{\mu}'(u)^2-2h_0(\mu)V_{\mu}''(u)=z^{-2F}V_1(z,\mu)^2-2h_0(\mu)V_2(z,\mu)$, with $z=\phi^{-1}(u)$ and where~$V_1,V_2\in\R[z,\mu]$. Clearly it suffices to show that 
 \begin{equation*}
  L(z,\mu)\!:=z^{-2F}V_1(z,\mu)^2-2h_0(\mu)V_2(z,\mu)> 0\text{ at $z=\phi^{-1}(u_r(\mu))=1-p_1(\mu)=1+\frac{b+\sqrt{b^2-4ac}}{2a}$}
 \end{equation*}
for all $\mu\in\Lambda$. 
(Here recall that $a,b,c\in\R(\mu)$ are the coefficients of the quadratic polynomial $q_{\mu}$ in the first integral~\refc{first-integral} and that $h_0\in R(\mu)$ is the energy of the outer boundary.) 

We claim that if $L(1-p_1(\hat\mu),\hat\mu)=0$ for some $\hat\mu=(\hat D,\hat F)\in\Lambda$ then the derivative of $D\longmapsto L(1-p_1(\mu),\mu)$ is strictly negative at $\mu=\hat\mu.$ To show this note first that $L(1-p_1(\hat\mu),\hat\mu)=0$ implies
 \begin{equation}\label{eq4}
  (1-p_1(\hat\mu))^{-2F}=\frac{2h_0(\hat\mu)V_2(1-p_1(\hat\mu),\hat\mu)}{V_1(1-p_1(\hat\mu),\hat\mu)^2}.
 \end{equation}
Moreover the derivative of $L(1-p_1(\mu),\mu)$ with respect $D$ is
 \begin{align*}
  \frac{d}{dD}L(1-p_1(\mu),\mu)=&\left.z^{-2F}\left(2V_1(z,\mu)\partial_zV_1(z,\mu)-\frac{2FV_1(z,\mu)^2}{z}\right)\right|_{z=1-p_1(\mu)}\partial_D \bigl(1-p_1(\mu)\bigr) \\[4pt] 
  &+\left.z^{-2F}\partial_D\left(V_1(z,\mu)^2\right)-2\partial_D\bigl(h_0(\mu)V_2(z,\mu)\bigr)\right|_{z=1-p_1(\mu)}.
 \end{align*}
The substitution of \refc{eq4} in the above expression evaluated at $\mu=\hat\mu$ gives us an expression of $\frac{d}{dD}L(1-p_1(\hat\mu),\hat\mu)$ which is algebraic in $\hat D$ and $\hat F.$ This is the key point in the proof. In doing so, with the help of an algebraic manipulator one can verify that
 \[
  \frac{d}{dD}L(1-p_1(\hat\mu),\hat\mu)=\frac{r_1(\hat\mu)\sqrt{\eta(\hat\mu)}+r_2(\hat\mu)}{r_3(\hat\mu)\sqrt{\eta(\hat\mu)}+r_4(\hat\mu)}\text{ with $\eta(\mu)\!:=\bigl(b^2-4ac\bigr)(\mu)$ and some $r_i\in\R[\mu].$}
  \]
Moreover
 \begin{align*}
  &r_1(\hat\mu)^2\eta(\hat\mu)-r_2(\hat\mu)^2=(1-2\hat F)^7(\hat D-\hat F+1)\hat D^4(\hat D+1)^3(\hat F+\hat D)(5\hat F\hat D-3\hat F^2+3\hat F-\hat D)\neq 0
  \intertext{and}
  &r_3(\hat\mu)^2\eta(\hat\mu)-r_4(\hat\mu)^2=-4\hat D^6\hat F^3(2\hat F-1)^9(\hat F+1)^3(\hat D+1)^3\neq 0
 \end{align*}
for all $(\hat D,\hat F)\in\Lambda.$ This proves that $\frac{d}{dD}L(1-p_1(\hat\mu),\hat\mu)$ is well defined and non-vanishing for all $\hat\mu\in\Lambda.$ Finally, since $\frac{d}{dD}L(1-p_1(\hat\mu),\hat\mu)<0$ at $\hat\mu=(-0.6,1.3),$ the claim follows. 

The claim implies that, for each fixed $F\in (1,3/2)$, the map $D\longmapsto L(1-p_1(\mu),\mu)$ has at most one zero for $D\in (-F,-1/2)$. This fact, on account of 
\[
 \left.L(1-p_1(\mu),\mu)\right|_{\mu=(-F,F)}=\frac{1}{F^2}>0\text{ and }
 \left.L(1-p_1(\mu),\mu)\right|_{\mu=(-\frac{1}{2},F)}=\frac{F^F(F-1)^{1-F}+2-3F}{4F^2(F-1)}>0
 \]
for all $F\in (1,3/2)$, shows the validity of the result.
\end{prova}

The following result gives the quantifier of the function $u\longmapsto (h_0(\mu)-V_{\mu}(u))\mathscr R_{\mu}(u)$, where
 \begin{equation}\label{eq11}
 \mathscr R_{\mu}\!:=\frac{(V_{\mu}')^2-2V_{\mu}V_{\mu}''}{(V_{\mu}')^3},
 \end{equation}
at the endpoints of the interval $\mathcal I_{\mu}=\left(-\frac{1}{F},u_r(\mu)\right).$

\begin{lema}\label{lema:quantificar_R}
The following holds:
\begin{enumerate}[$(a)$]
\item If $\hat\mu\in\Lambda\setminus\{F=2\}$ then 
$(h_0(\mu)-V_{\mu}(u))\mathscr R_{\mu}(u)\sim_{-\frac{1}{F}} \frac{(F-2)h_0(\mu)}{D(F-1)}(Fu+1)^{\frac{2}{F}-1}$ at $\hat\mu$.

\item If $\hat\mu\in\Lambda$ then $(h_0(\mu)-V_{\mu}(u))\mathscr R_{\mu}(u)\sim_{u_r(\mu)}\mathscr R_{\mu}(u_r(\mu))V'_{\mu}(u_r(\mu))(u_r(\mu)-u)$ at $\hat\mu$.

\end{enumerate}

\end{lema}

\begin{prova}
To prove $(a)$ let us fix $\hat\mu\in\Lambda\setminus\{F=2\}$. Taking \eqref{eq:expressions_V} into account, and with the help of an algebraic manipulator, one can verify that
 \[
(h_0(\mu)-V_{\mu}(u))\mathscr R_{\mu}(u)=\left.\frac{V_0(z,\mu)z^{-F}f_{\mu}(z)}{V_1^3(z,\mu)}\right|_{z=\phi^{-1}(u)},
\]
where $f_{\mu}(z)\!:=V_1(z,\mu)^2+2V_2(z,\mu)\bigl(V_0(z,\mu)-h_0(\mu)z^{2F}\bigr)$ is the sum of $7$ monomials of the form $c(\mu)z^{n_1+n_2F}$ with $n_i\in\Z$ for $i=1,2$, and $c(\mu)$ a well defined rational function at $\mu=\hat\mu$. In addition, the monomial with the largest exponent for $\mu\approx\hat\mu$ is $\frac{D(1+D-F)(F-2)}{F(F-1)(2F-1)}z^{2+2F}$. Accordingly 
\[
f_{\mu}(z)\sim_{+\infty} \frac{D(1+D-F)(F-2)}{F(F-1)(2F-1)}z^{2+2F} \text{ at }\hat\mu.
\]
On the other hand, taking \eqref{eq:expressions_V} into account once again,
$V_0(z,\mu)\sim_{+\infty} \frac{D}{2-2F}z^2$ and $V_1(z,\mu)\sim_{+\infty} Dz^2$ at $\hat\mu.$
Consequently
\[
\frac{V_0(z,\mu)z^{-F}f_{\mu}(z)}{V_1^3(z,\mu)}\sim_{+\infty} \frac{(F-D-1)(F-2)}{2DF(F-1)^2(2F-1)} z^{F-2} \text{ at }\hat\mu,
\] 
and by applying Lemma~\ref{lema:traduccio} we get that
\[
(h_0(\mu)-V_{\mu})\mathscr R_{\mu}(u)\sim_{-1/F}
\frac{(F-D-1)(F-2)}{2DF(F-1)^2(2F-1)} (Fu+1)^{\frac{2}{F}-1} \text{ at }\hat\mu.
\]
Since $h_0(\mu)=\frac{F-D-1}{2F(F-1)(2F-1)}$, this proves $(a)$. To show $(b)$ note that $V_{\mu}'(u)^2-2V_{\mu}(u)V_{\mu}''(u)$ does not vanish
at $u=u_r(\mu)$ by \lemc{lema:corba}, and that $h_0(\mu)-V_{\mu}(u) \sim_{u_r(\mu)} V_{\mu}'(u_r(\mu))(u_r(\mu)-u)$ by $(b)$ in \lemc{lema:quantificar_V}. This proves the result. 
\end{prova}

In order to state our next result let us define
\[
\Psi_{\mu}(u)\!:=\frac{1}{V_{\mu}'(u)}W\!\!\left[\left(\frac{V_{\mu}}{h_0(\mu)-V_{\mu}}\right)^{\frac{1}{2}\nu(\mu)},(h_0(\mu)-V_{\mu})V_{\mu}^{\frac{1}{2}}\mathscr R_{\mu}\right]\!\!(u),
\]
where $\mathscr R_{\mu}$ is the function in \refc{eq11} and $\nu:\Lambda\rightarrow\R$ is a continuous function to be determined. Note that $\Psi_{\mu}$ is the function in $(ii)$ of \teoc{teorema-jdde} for the particular case $n=1$. Let us advance that $\nu$ is to be chosen in such a way that the family \eqref{familia} is continuously quantifiable at $z=1$, so that we can apply \teoc{teorema-jdde}.

\begin{prop}\label{prop:quantificar_Psi}
Let us fix $\hat\mu=(\hat D,\hat F)\in\Lambda.$ Then the following holds:
\begin{enumerate}[$(a)$]
\item If $\hat F\neq 2$ and $\nu(\hat\mu)\neq \frac{\hat F-2}{\hat F-1}$ then $\Psi_{\mu}(u)\sim_{-\frac{1}{F}} a(\mu)(Fu+1)^{\frac{4+\nu-F(3+\nu)}{F}}$ at $\hat\mu$, where
\[
 a(\mu)\!:=-(F-2)(F-2-\nu(F-1))(-D)^{\frac{\nu}{2}}(F-1)^{\nu-\frac{5}{2}}F^{\frac{\nu+1}{2}}(F-D-1)^{\frac{\nu+3}{2}}(4F-2)^{\frac{\nu-3}{2}}.
\]
\item If $\nu(\hat\mu)\neq -2$ then $\Psi_{\mu}(u)\sim_{u_r(\mu)} b(\mu)(u_r(\mu)-u)^{-\frac{\nu}{2}}$ at $\hat\mu$, where
$b(\mu)\!:=-\frac{\nu+2}{2}\mathscr R_{\mu}(u_r)h_0^{\frac{\nu+1}{2}}V_{\mu}'(u_r)^{-\frac{\nu}{2}}.$
\end{enumerate}
\end{prop}

\begin{prova}
Let us fix $\hat\mu=(\hat D,\hat F)\in\Lambda$. A computation shows that
\begin{equation}\label{eq:Psi}
\Psi_{\mu}=\frac{\psi_{\mu}}{2V_{\mu}^{\frac{1}{2}}(V_{\mu}')^5}\left(\frac{V_{\mu}}{h_0(\mu)-V_{\mu}}\right)^{\frac{\nu(\mu)}{2}},
\end{equation}
where, omitting the dependence on $\mu$ for the sake of shortness,
\begin{equation*}
\psi_{\mu}\!:=4V^2(V-h_0)V'V'''-(V'^2-2VV'')\left(V'^2(h_0(\nu-1)+3V)+6(h_0-V)VV''\right).
\end{equation*}
Taking \refc{eq:expressions_V} into account, and with the help of an algebraic manipulator, one can show that $\psi_{\mu}(u)=f_{\mu}(z)$ where $z=\phi^{-1}(u)=(Fu+1)^{-1/F}$ and $f_{\mu}$ is the sum of~25 monomials of the form $c(\mu)z^{n_1+n_2F}$ with $n_i\in\Z$ and $c(\mu)$ a well defined rational function at~$\mu=\hat\mu$. In particular, the monomial with largest exponent for $\mu\approx\hat\mu$ is $c(\mu)z^{6-2F}$ with
 \[
 c(\mu)\!:=-\frac{(F-2)(F-2-\nu(F-1))D^3(1+D-F)^2}{2F^2(F-1)^3(1-2F)^2}.
 \]
Accordingly, if $\hat F\neq 2$ and we choose any $\nu$ such that $\nu(\hat\mu)\neq \frac{\hat F-2}{\hat F-1}$, we get that 
$f_{\mu}(z)\sim_{+\infty}c(\mu)z^{6-2F}$ at~$\hat\mu.$ Thus, due to $\psi_{\mu}(u)=f_{\mu}(z)$ with $z=\phi^{-1}(u)=(Fu+1)^{-1/F}$, by applying \lemc{lema:traduccio} we can assert that 
 \begin{equation*}
  \psi_{\mu}(u)\sim_{-\frac{1}{F}} c(\mu)(Fu+1)^{2-6/F}\text{ at $\hat\mu.$}
 \end{equation*} 
On the other hand, from \eqref{eq:expressions_V} and applying \lemc{lema:traduccio} again, we get
$V_{\mu}(u)\sim_{-\frac{1}{F}}h_0(\mu)=\frac{F-D-1}{2F(F-1)(2F-1)}$ and $V_{\mu}'(u)\sim_{-\frac{1}{F}}D(Fu+1)^{1-2/F}$ at $\hat\mu.$ Taking these three quantifiers into account, together with the quantifier for $h_0(\mu)-V_{\mu}(u)$ given by $(a)$ in \lemc{lema:quantificar_V}, from \refc{eq:Psi} we can assert that $\Psi_{\mu}(u)\sim_{-\frac{1}{F}} a(\mu)(Fu+1)^{\frac{4+\nu-F(3+\nu)}{F}}$ at $\hat\mu$ with
\[
 a(\mu)=-(F-2)(F-2-\nu(F-1))(-D)^{\frac{\nu}{2}}(F-1)^{\nu-\frac{5}{2}}F^{\frac{\nu+1}{2}}(F-D-1)^{\frac{\nu+3}{2}}(4F-2)^{\frac{\nu-3}{2}}.
\]
This shows $(a)$. Let us turn now to the proof of the claim in $(b)$. Since $\lim_{(\mu,u)\to (\hat\mu,u_r(\hat\mu))}V_{\mu}(u)=h_0(\hat\mu)$, by using \lemc{lema:corba} and $(b)$ in \lemc{lema:quantificar_V} we obtain that
\[
\psi_{\mu}(u)\sim_{u_r(\mu)} (\nu+2)\left(2V_{\mu}(u_r)V_{\mu}''(u_r)-V_{\mu}'(u_r)^2\right)V_{\mu}'(u_r)^2h_0 \text{ at }\hat\mu,
\]
provided that $\nu(\hat\mu)\neq -2$. On account of this, exactly the same ingredients yield to
\[
\Psi_{\mu}(u)\sim_{u_r(\mu)} -\frac{\nu+2}{2}\mathscr R_{\mu}(u_r)h_0^{\frac{\nu+1}{2}}V_{\mu}'(u_r)^{-\frac{\nu}{2}}(u_r-u)^{-\frac{\nu}{2}} \text{ at }\hat\mu,
\]
as we desired. This concludes the proof of the result. 
\end{prova}
%
%
%

\begin{prooftext}{Proof of \teoc{thm:crit}.}
We begin by showing that any $\hat\mu=(\hat D,\hat F)\in\Lambda\setminus(\Gamma_B\cup\Gamma_U)$ is a local regular value of the period function at the outer boundary. This will follow by applying \teoc{teorema-jdde} with $n=0$. With this aim in view let $f_{\mu}(z)$ be twice the even part of the function 
 \begin{equation}\label{eq17}
  z\longmapsto z\sqrt{h_0(\mu)}\left(g_{\mu}^{-1}\right)''\!\bigl(z\sqrt{h_0(\mu)}\,\bigr).
 \end{equation} 
Note that $f_{\mu}(z)$ is precisely the function in \refc{familia} for $n=0$ because $\DD_{\boldsymbol\nu_0}=id$ by definition. That being said, we first apply the second part of \teoc{teorema-jdde}, which requires the quantifiers of $h_0(\mu)-V_{\mu}$ and $(h_0(\mu)-V_{\mu})V_{\mu}^{1/2}\mathscr R_{\mu}$ at the endpoints of $\left(-\frac{1}{F},u_r(\mu)\right)$. In this regard \lemc{lema:quantificar_V} shows that the family $\{h_0(\mu)-V_{\mu}\}_{\mu\in\Lambda}$ is continuously quantifiable in $\hat\mu$ at $u=-\frac{1}{F}$ by $\beta_{\ell}(\mu)=\frac{2-2F}{F}$ with limit $b_{\ell}(\mu)=\frac{DF^{2-\frac{2}{F}}}{2(1-F)}$ and at $u=u_r(\mu)$ by $\beta_r(\mu)=-1$ with limit $b_r(\mu)=V_{\mu}'(u_r(\mu))$. On the other hand, from \lemc{lema:quantificar_R} it follows that the family $\{(h_0(\mu)-V_{\mu})V_{\mu}^{1/2}\mathscr R_{\mu}\}_{\mu\in\Lambda}$ is continuously quantifiable in $\hat\mu$ at $u=-\frac{1}{F}$ by $\alpha_{\ell}(\mu)=1-\frac{2}{F}$ with limit $a_{\ell}(\mu)=\frac{(F-2)h_0(\mu)^{\frac{3}{2}}F^{\frac{2}{F}-1}}{D(F-1)}$ and at $u=u_r(\mu)$ by $\alpha_r(\mu)=-1$ with limit $a_r(\mu)=\mathscr R_{\mu}(u_r(\mu))V'_{\mu}(u_r(\mu))h_0(\mu)^{\frac{1}{2}}$. Thus $\frac{\alpha_{\ell}}{\beta_{\ell}}(\mu)=\frac{2-F}{2(F-1)}$ and $\frac{\alpha_{r}}{\beta_{r}}(\mu)=1,$ and by applying the second part of \teoc{teorema-jdde} with $n=0$ we get that $\{f_{\mu}\}_{\mu\in\Lambda}$ is continuously quantifiable at $z=1$ by 
\[
\xi(\mu)=-\min\left\{\frac{2-F}{2(F-1)},1\right\}+1=
\begin{cases}
0 & \text{ if }1<F\leq \frac{4}{3},\\
\frac{3F-4}{2(F-1)} & \text{ if } F>\frac{4}{3}.
\end{cases}
\]
Here it is to be pointed out that in order to cover the case $F=\frac{4}{3}$ we checked that $\frac{a_{r}}{b_{r}}(\mu)+\frac{a_{\ell}}{b_{\ell}}(\mu)\neq 0$ for all $\mu=(D,\frac{4}{3}).$ Indeed, note that $\frac{a_{r}}{b_{r}}(\mu)=\mathscr R_{\mu}(u_r(\mu))h_0(\mu)^{\frac{1}{2}},$ which is strictly positive by \lemc{lema:corba} and the definition of $\mathscr R_{\mu}$ in~\eqref{eq11}. On the other hand, if $F=\frac{4}{3}$ then a computation shows that $\frac{a_{\ell}}{b_{\ell}}(\mu)=\frac{4h_0(\mu)^{\frac{3}{2}}}{3D^2}>0$. Accordingly we have $\frac{a_{r}}{b_{r}}(\mu)+\frac{a_{\ell}}{b_{\ell}}(\mu)> 0$ for all $\mu=(D,\frac{4}{3})$.

If $\hat F>\frac{3}{2}$ then $\xi(\hat\mu)>1/2$, and by applying $(a)$ in \teoc{teorema-jdde} we have $\mathrm{Crit}\bigl((\out_{\hat\mu},X_{\hat\mu}),X_{\mu}\bigr)=0$. Note also that in this case, as it is explained in \obsc{remark}, the sign of the derivative of the period function near the outer boundary is given by $a_{\ell}(\mu)=\frac{(F-2)h_0(\mu)^{\frac{3}{2}}F^{\frac{2}{F}-1}}{D(F-1)},$ that changes at $F=2.$ Using Bolzano's Theorem, this easily implies that $\mathrm{Crit}\bigl((\out_{\hat\mu},X_{\hat\mu}),X_{\mu}\bigr)\geq 1$ for all $\hat\mu=(\hat D,2)\in\Lambda.$

If $\hat F\in (1,\frac{3}{2})$ then $0\leq \xi(\hat\mu)<\frac{1}{2}$, and so we need to apply $(b)$ in \teoc{teorema-jdde}, which lead us to the computation of $N_1(\mu)$, the first momentum of $f_{\mu}$. For the sake of shortness, to this end we take advantage of the results in \cite{ManRojVil2016-JDE,MMV2}. (Let us remark however that, although is a lengthly computation, it could be done without appealing to these results.) In doing so we get
 \begin{align*}
  N_1(\mu)=&\int_0^1\frac{f_{\mu}(z)}{\sqrt{1-z^2}}dz
  =\sqrt{h_0(\mu)}\int_{-1}^1\frac{z\!\left(g_{\mu}^{-1}\right)''\!\bigl(\sqrt{h_0(\mu)}z\bigr)}{\sqrt{1-z^2}}dz \\[3pt]
  =&\sqrt{h_0(\mu)}\int_{-\frac{\pi}{2}}^{\frac{\pi}{2}}\!\left(g_{\mu}^{-1}\right)''\!\bigl(\sqrt{h_0(\mu)}\sin\theta\bigr)\sin\theta d\theta
  =\lim_{h\rightarrow h_0(\mu)}\frac{T'_{\mu}(h)}{\sqrt{2}h_0(\mu)},
 \end{align*}
where the first equality follows from \defic{def_moment}, the second one using that $f_{\mu}(z)$ is twice the even part of the function in~\refc{eq17} and the last one by Corollary 3.12 in \cite{ManRojVil2016-JDE}. Note also that the first integral above is convergent thanks to $\xi(\mu)<1/2$. We take now advantage of \propc{prop:corba_bif}, which shows that if $F\in (1,\frac{3}{2})$ then $\lim_{s\to 0^+}P'(s;\mu)=\Delta(\mu)$. Recall that $P(s;\mu)$ refers to the period of the periodic orbit of $\refc{Loud2}$ passing through the point $\bigl(p_1(\mu)-s,0\bigr)$ and so it is clear that $P(s;\mu)=T_{\mu}(\zeta(s;\mu))$, where $s\longmapsto\zeta(s;\mu)$ is an orientation reversing diffeomorphism. (Actually $\zeta(s;\mu)=V_{\mu}\bigl(\phi(1-p_1+s)\bigr)$ but this is not relevant for our purposes.) Hence $N_1(\mu)=-C(\mu)\Delta(\mu)$, where $C(\mu)$ is positive and $\Delta(\mu)$ is the coefficient given in \propc{prop:corba_bif}. In particular $N_1(\mu)=0$ if and only if $D=\mathcal G(F)$. Thus, if $\hat\mu\in\Lambda\setminus\Gamma_B$ then $N_1(\hat\mu)\neq 0$ and so, by applying $(b1)$ in \teoc{teorema-jdde} with $j=1,$ we get that $\mathrm{Crit}\bigl((\out_{\hat\mu},X_{\hat\mu}),X_{\mu}\bigr)=0$ as desired. Exactly as before, the fact that $\text{Crit}((\Pi_{\hat\mu},X_{\hat\mu}),X_{\mu})\geq 1$ for any $\hat\mu\in\Gamma_B$ is because $N_1(\mu)$, and then the derivative of the period function, changes sign as we cross the curve $D=\mathcal G(D).$ 

It only remains to be proved that $\text{Crit}((\Pi_{\hat\mu},X_{\hat\mu}),X_{\mu})\leq 1$ for any $\hat\mu\in\Gamma_B$ with $\hat F\in (\frac{4}{3},\frac{3}{2})$. To this end we shall apply \teoc{teorema-jdde} with $n=1$, and so we need first to choose a convenient function $\nu_1(\mu)$ in such a way that the family $\{\DD_{\nu_1(\mu)}\!\circ\! f_{\mu}\}_{\mu\in\Lambda}$ 
is continuously quantifiable at $z=1.$ With this aim in view recall that $\{f_{\mu}\}_{\mu\in\Lambda}$ is continuously quantifiable at $z=1$ by $\xi(\mu)<\frac{1}{2}$ for the parameter values under consideration. Therefore we can apply \lemc{lem:reduc_moment_finit} with $n=\ell=1$, which shows that $N_1\left[\mathscr D_{\nu_{1}}\!\circ\! f_{\mu}\right]=-(\nu_1+1)N_1[f_{\mu}]$ for any $\nu_1\in\R$. Accordingly, if we define $\nu_1(\mu)=-1$ for all $\mu$ then 
the first momentum of $\mathscr D_{\nu_{1}}\!\circ\! f_{\mu}$ is identically zero. With this choice for $\nu_1(\mu)$ we turn to the quantification of $\{\DD_{\nu_1(\mu)}\!\circ\! f_{\mu}\}_{\mu\in\Lambda}$ at $z=1$ by means of the second part of \teoc{teorema-jdde}. To do so we need the quantifiers of $h_0(\mu)-V_{\mu}$ and 
\[
u\longmapsto \frac{1}{V_{\mu}'(u)}W\!\!\left[\left(\frac{V_{\mu}}{h_0(\mu)-V_{\mu}}\right)^{-\frac{1}{2}},(h_0(\mu)-V_{\mu})V_{\mu}^{\frac{1}{2}}\mathscr R_{\mu}\right]\!\!(u),
\]
at the endpoints of $\left(-\frac{1}{F},u_r(\mu)\right)$. To obtain the ones of the later we apply \propc{prop:quantificar_Psi} with $\nu=-1$, which shows that it is continuously quantifiable at $u=-\frac{1}{F}$ by $\alpha_{\ell}(\mu)=\frac{2F-3}{F}$ and at $u=u_r(\mu)$ by $\alpha_r(\mu)=-\frac{1}{2}$. As we already used, \lemc{lema:quantificar_V} shows that the family $\{h_0(\mu)-V_{\mu}\}_{\mu\in\Lambda}$ is continuously quantifiable at $u=-\frac{1}{F}$ by $\beta_{\ell}(\mu)=\frac{2-2F}{F}$  and at $u=u_r(\mu)$ by $\beta_r(\mu)=-1$. In this case $\frac{\alpha_{\ell}}{\beta_{\ell}}(\mu)=\frac{2F-3}{2(1-F)}$ and $\frac{\alpha_{r}}{\beta_{r}}(\mu)=\frac{1}{2},$ and by applying the second part of \teoc{teorema-jdde} with $n=1$ we can assert that $\{\DD_{\nu_1(\mu)}\!\circ\! f_{\mu}\}_{\mu\in\Lambda}$ is continuously quantifiable at $z=1$ by 
\[
\xi(\mu)=-\min\left\{\frac{2F-3}{2(1-F)},\frac{1}{2}\right\}+\frac{1}{2}=
\begin{cases}
\frac{3F-4}{2(F-1)}& \text{ if }F>\frac{4}{3},\\
0 & \text{ if }F\in(1,\frac{4}{3}),
\end{cases}
\] 
In particular, $\xi(\mu)\in (0,\frac{1}{2})$ for all $F\in (\frac{4}{3},\frac{3}{2})$. Thus, by applying $(b2)$ in \teoc{teorema-jdde} with $n=m=1$, due to $N_1\left[\mathscr D_{\nu_{1}}\!\circ\! f_{\mu}\right]\equiv 0$, we can assert that $\text{Crit}\bigl((\Pi_{\hat\mu},X_{\hat\mu}),X_{\mu}\bigr)\leq 1$ as desired.
\end{prooftext}

Let us conclude the paper by making some comments regarding the difficulties we have encountered in the study of the criticality of other parameters in $\Gamma_B$ apart from the ones contemplated in \teoc{thm:crit}. As we already mentioned, we apply the tools developed in~\cite{ManRojVil2016-JDE,ManRojVil2016-JDDE} to bound the criticality of a polycycle $\Pi_{\hat\mu}$ in a family of potential systems $Y_{\mu}=-y\partial_x+V'_{\mu}(x)\partial_y$. In short, they apply to two different settings: either the energy $h_0(\mu)$ of the potential function $V_{\mu}$ at the polycycle is finite for all $\mu\approx\hat\mu$ or  $h_0(\mu)=+\infty$ for all $\mu\approx\hat\mu.$ (For each one of these situations we have a specific result with its own hypothesis to be verified, see respectively Theorems B and A in \cite{ManRojVil2016-JDDE}.) We cannot treat the case in which in any neighbourhood of $\hat\mu$ there are $\mu_1$ and $\mu_2$ with $h_0(\mu_1)<+\infty$ and $h_0(\mu_2)=+\infty.$ This is precisely what happens in the segments $\{\mu\in\Gamma_B:(F+D)D=0\}$. On the contrary we can apply the mentioned theorems to the rest of the parameters in $\Gamma_B$, but the technical hypothesis are not verified for different reasons. Indeed, for $\{F=\frac{1}{2}$\} the first integral~\refc{first-integral} has a pole and this makes $\{h_0(\mu)-V_{\mu}\}_{\mu}$ not quantifiable. The application of \teoc{teorema-jdde} with $n=1,2$ for $\hat\mu\in\{F=2\}$ yields to $\xi(\hat\mu)=\frac{1}{2}.$ Finally, if $\hat\mu\in\{D=-\frac{1}{2}\}$ then $\xi(\hat\mu)<1/2$, $m=1$ and $N_1\equiv 0$, which is far from being understood and poses additional difficulties.


\begin{thebibliography}{99}

\bibitem{ChiJac} C. Chicone and M. Jacobs, \emph{Bifurcation of critical periods for plane vector fields,} Trans. Amer. Math. Soc. {\bf 312} (1989) 433--486. 

\bibitem{Chicone} C.~Chicone,  \emph{Review in MathSciNet, ref. 94h:58072.}

\bibitem{Chouikha} A.R. Chouikha, \emph{Monotonicity of the period function for some planar differential systems. I. Conservative and quadratic systems,} Appl. Math. (Warsaw) {\bf 32} (2005) 305--325. 

\bibitem{CG} W.A.~Coppel and L.~Gavrilov,  \emph{The period function of a Hamiltonian quadratic system,}
Differential Integral Equations {\bf 6} (1993) 1357--1365.

\bibitem{GGV04} A. Gasull, A. Guillamon and J. Villadelprat, \emph{The period function for second-order quadratic ODEs is monotone,} Qual. Theory Dyn. Syst. {\bf 5} (2004) 201--224.

\bibitem{Ilyashenko} Y. Ilyashenko, \emph{Centennial history of Hilbert's 16th problem,} Bulletin of the AMS, {\bf 39} (2002) 301--354. 

\bibitem{Hsu} Sze-Bi Hsu, \emph{A remark on the period of the period solution in the Lotka-Volterra system}, 
J. Math. Anal. Appl. {\bf 95} (1983) 428--436.

\bibitem{ManRojVil2016-JDE} F.~Ma{\~n}osas, D.~Rojas and J.~Villadelprat, \emph{The criticality of centers of potential systems at the outer boundary,} J. Differential Equations, {\bf 260} (2016) 4918--4972.

\bibitem{ManRojVil2016-JDDE}
F.~Ma{\~n}osas, D.~Rojas and J.~Villadelprat,
\emph{Analytic tools to bound the criticality at the outer boundary of the period annulus}, J. Dynamics and Differential Equations (2016). doi:10.1007/s10884-016-9559-x

\bibitem{TopaV} F. Ma–osas and J. Villadelprat,
\emph{The bifurcation set of the period function of the dehomogenized Loud's centers is bounded,}
Proc. Amer. Math. Soc. {\bf 136} (2008) 1631--1642. 

\bibitem{MMV2}
P.~Marde{\v{s}}i{\'c}, D.~Mar{\'{\i}}n and J.~Villadelprat,
\emph{The period function of reversible quadratic centers}, J. Differential Equations {\bf 224} (2006) 120--171.

\bibitem{MMSV} P.~Marde{\v{s}}i{\'c}, D. Mar'n, M. Saavedra and J. Villadelprat, 
\emph{Unfoldings of saddle-nodes and their Dulac time},
J. Differential Equations {\bf 261} (2016) 6411--6436.

\bibitem{MV} D. Mar'n and J. Villadelprat,
\emph{On the return time function around monodromic polycycles,}
J. Differential Equations {\bf 228} (2006) 226--258. 

\bibitem{Rothe}
F.~Rothe, 
\emph{The periods of the Volterra-Lotka system,} 
J. Reine Angew. Math. {\bf 355} (1985) 129--138.

\bibitem{Roussarie} R. Roussarie, ``Bifurcation of Planar Vector Fields and HilbertÕs Sixteenth Problem'', Progress in
Mathematics, vol. 164, BirkhŠuser, Basel, 1998.

\bibitem{Jordi} J.~Villadelprat, \emph{On the reversible quadratic centers with monotonic period function,}
Proc. Amer. Math. Soc. {\bf 135} (2007) 2555--2565.

\bibitem{Waldvogel}
J.~Waldvogel, 
\emph{The period in the Lotka-Volterra system is monotonic,} 
J. Math. Anal. Appl. {\bf 114} (1986) 178--184.

\bibitem{Zhao} Y.~Zhao,  \emph{The monotonicity of period function for codimension four quadratic system,}
J. Differential Equations {\bf 185} (2002) 370--387.

\bibitem{Zhao2}  Y. Zhao, \emph{The period function for quadratic integrable systems with cubic orbits,}
J. Math. Anal. Appl. {\bf 301} (2005) 295--312. 

\bibitem{Zhao3}  Y. Zhao, \emph{On the monotonicity of the period function of a quadratic system,}
Discrete Contin. Dyn. Syst. {\bf 13} (2005) 795--810. 

\bibitem{Zhao5} H. Liang and Y. Zhao, \emph{On the period function of reversible quadratic centers
with their orbits inside quartics,} Nonlinear Anal. {\bf 71} (2009) 5655--5671. 

\bibitem{Zhao4} H. Liang and Y. Zhao, \emph{On the period function of a class of reversible quadratic centers,} 
Acta Math. Sin. (Engl. Ser.) {\bf 27} (2011) 905--918. 

\end{thebibliography}
\end{document}